\documentclass[a4paper,12pt]{amsart}
\usepackage{amssymb}
\usepackage{ifthen}
\usepackage{graphicx}
\usepackage{subfigure}
\usepackage{color,xcolor}
\nonstopmode
\numberwithin{equation}{section}
\newtheorem{thm}{Theorem}
\newtheorem{cor}{Corollary}
\newtheorem{lem}{Lemma}
\newtheorem{prop}{Proposition}

\newtheorem{conj}{Conjecture}
\newtheorem{prob}{Problem}

\theoremstyle{definition}
\newtheorem{defn}{Definition}
\newtheorem{ca}{Case}

\newtheorem{rem}{Remark}

\newenvironment{pf}[1][]{%
 \vskip 1mm
 \noindent
 \ifthenelse{\equal{#1}{}}%
  {{\slshape Proof. }}%
  {{\slshape #1.} }%
 }%
{\qed\medskip}

\newcounter{alphabet}

\newenvironment{Thm}[1][]{\refstepcounter{alphabet}%
\bigskip%
\noindent%
{\bf Theorem \Alph{alphabet}}%
{\bf .} \itshape}{\vskip 8pt}

\newcommand{\ID}{{\mathbb D}}

\def\be{\begin{equation}}
\def\ee{\end{equation}}

\newcommand{\ben}{\begin{enumerate}}
\newcommand{\een}{\end{enumerate}}

\newcommand{\blem}{\begin{lem}}
\newcommand{\elem}{\end{lem}}
\newcommand{\bthm}{\begin{thm}}
\newcommand{\ethm}{\end{thm}}
\newcommand{\bcor}{\begin{cor}}
\newcommand{\ecor}{\end{cor}}
\newcommand{\beg}{\begin{exam}}
\newcommand{\eeg}{\end{exam}}
\newcommand{\begs}{\begin{examples}}
\newcommand{\eegs}{\end{examples}}
\newcommand{\bdefe}{\begin{defn}}
\newcommand{\edefe}{\end{defn}}
\newcommand{\bprob}{\begin{prob}}
\newcommand{\eprob}{\end{prob}}
\newcommand{\bques}{\begin{ques}}
\newcommand{\eques}{\end{ques}}
\newcommand{\bei}{\begin{itemize}}
\newcommand{\eei}{\end{itemize}}
\newcommand{\bcon}{\begin{conj}}
\newcommand{\econ}{\end{conj}}
\newcommand{\bop}{\begin{op}}
\newcommand{\eop}{\end{op}}

\newcommand{\bas}{\begin{assertion}}
\newcommand{\eas}{\end{assertion}}

\newcommand{\bfa}{\begin{fact}}
\newcommand{\efa}{\end{fact}}

\newcommand{\bca}{\begin{ca}}
\newcommand{\eca}{\end{ca}}

\newcommand{\bst}{\begin{step}}
\newcommand{\est}{\end{step}}

\newcommand{\bsca}{\begin{sca}}
\newcommand{\esca}{\end{sca}}

\newcommand{\bcl}{\begin{cl}}
\newcommand{\ecl}{\end{cl}}

\newcommand{\bmlem}{\begin{mlem}}
\newcommand{\emlem}{\end{mlem}}

\newcommand{\bscl}{\begin{scl}}
\newcommand{\escl}{\end{scl}}

\newcommand{\bcons}{\begin{conjs}}
\newcommand{\econs}{\end{conjs}}

\newcommand{\bprop}{\begin{prop}}
\newcommand{\eprop}{\end{prop}}

\newcommand{\br}{\begin{rem}}
\newcommand{\er}{\end{rem}}
\newcommand{\brs}{\begin{rems}}
\newcommand{\ers}{\end{rems}}
\newcommand{\bo}{\begin{obser}}
\newcommand{\eo}{\end{obser}}
\newcommand{\bos}{\begin{obsers}}
\newcommand{\eos}{\end{obsers}}
\newcommand{\bpf}{\begin{pf}}
\newcommand{\epf}{\end{pf}}
\newcommand{\ba}{\begin{array}}
\newcommand{\ea}{\end{array}}
\newcommand{\beq}{\begin{eqnarray}}
\newcommand{\beqq}{\begin{eqnarray*}}
\newcommand{\eeq}{\end{eqnarray}}
\newcommand{\eeqq}{\end{eqnarray*}}

\newcommand{\ds}{\displaystyle}

\newcounter{minutes}
\divide\time by 60
\newcounter{hours}
\multiply\time by 60 \addtocounter{minutes}{-\time}

\begin{document}

\bibliographystyle{amsplain}
\title [Landau-type theorems for certain bounded bi-analytic functions and biharmonic mappings]
{Landau-type theorems for certain bounded bi-analytic functions and biharmonic mappings}

\def\thefootnote{}
\footnotetext{ \texttt{\tiny File:~\jobname .tex,
          printed: \number\day-\number\month-\number\year,
          \thehours.\ifnum\theminutes<10{0}\fi\theminutes}
} \makeatletter\def\thefootnote{\@arabic\c@footnote}\makeatother

\author{Ming-Sheng Liu 
}
\address{M.S. Liu, School of Mathematical Sciences, South China Normal University, Guangzhou, Guangdong 510631, China.}
\email{liumsh65@163.com}

\author{Saminathan Ponnusamy}
\address{S. Ponnusamy, Department of Mathematics,
Indian Institute of Technology Madras, Chennai-600 036, India.}
\address{Lomonosov Moscow State University, Moscow Center of Fundamental and Applied Mathematics, Moscow, Russia. }
\email{samy@iitm.ac.in}



\begin{abstract}
In this paper,  we establish three new versions of Landau-type theorems for bounded bi-analytic functions of the form
$F(z)=\bar{z}G(z)+H(z)$, where $G$ and $H$ are analytic in the unit disk $|z|<1$ with $G(0)=H(0)=0$ and $H'(0)=1$. In particular,
two of them are sharp while the other one either generalizes or improves the corresponding result of Abdulhadi and Hajj.
As consequences, several new sharp versions of Landau-type theorems for certain subclasses of bounded biharmonic mappings are proved.
\end{abstract}

\maketitle
\pagestyle{myheadings}
\markboth{M.S. Liu and S. Ponnusamy}{Landau-type theorems for bi-analytic functions and biharmonic mappings}

\section{Introduction and Preliminaries}\label{HLP-sec1}

One of the open problems in classical complex analysis is to obtain the precise value of the Bloch constant for analytic
functions in the unit disk. In \cite{CGH2000}, Chen  et al, and Guathier considered the analogous problem of estimating the Bloch constant
for planar harmonic mappings. See also the work of Chen and Guathier \cite{cg} for  planar harmonic and pluriharmonic mappings.
Motivated by the work from \cite{CGH2000}, this topic was dealt by a number of authors with considerable improvements over the previously known Landau-type theorems. These will be indicated later in this section. In this article, we consider bi-analytic and biharmonic mappings and establish several new sharp versions of Landau-type theorems for these two classes of mappings.

\subsection{Definitions and Notations}
A complex-valued function $f$ is a bi-analytic (resp. harmonic) on a domain $D\subset\mathbb{C}$ if and only if $f$ is twice continuously differentiable and satisfies the bi-analytic equation
$f_{\bar{z}\bar{z}}(z)=0$ (resp. Laplacian equation $f_{z\overline{z}}(z)=0$) in  $D$, where we use the common notations for its formal derivatives:
$$
f_z=\frac{1}{2}(f_x - i f_y),~\mbox{ and }~f_{\overline{z}}=\frac{1}{2}(f_x + i f_y), \quad  z=x+iy.
$$
Note also that
$$ \Delta f=4f_{z\overline{z}}=\frac{\partial^2 f}{\partial x^2}+\frac{\partial^2 f}{\partial y^2}.
$$
It is well-known that every bi-analytic function $f$ in a simply connected domain $D$  has the representation (cf. \cite{AH2022})
$$
f(z)=\bar{z}g(z)+h(z),
$$
where $g$ and $h$ are complex-valued analytic functions in $D$.
%
Similarly, every harmonic function $f$ in a simply connected domain $D$  can be written as $f=h+\overline{g}$ with $f(0)=h(0)$,
where $g$ and $h$ are analytic on $D$ (for details see \cite{CS1984}).

A complex-valued function $F$ is said to be {\it biharmonic} on a domain $D\subset\mathbb{C}$ if and only if $F$ is four times continuously differentiable and satisfies  the biharmonic equation $\Delta (\Delta f)=0$ in $D$. It is well-known (cf. \cite{AMK2005}) that a   biharmonic mapping $F$ in a simply connected domain $D$ has the following representation:
$$
F(z)=|z|^2G(z)+H(z),
$$
where $G$ and $H$ are harmonic in $D$.


A domain $D\subset\mathbb{C}$ is said to be starlike  if and only if the line segment $[0,w]$ joining the origin $0$ to every other point $w\in D$ lies entirely in $D$.

\bdefe (cf. \cite{LP2021,M1980,PQ2013})
A continuously differentiable function $F$ on $D$ is said to be fully starlike in $D$ if it is sense-preserving, $F(0)= 0, F(z)\neq 0$ in $D\backslash\{0\}$ and the curve $F(re^{it})$ is starlike with respect to the origin for each $r\in (0, 1)$. The last condition is
same as saying that
\begin{eqnarray*}
 \frac{\partial\arg F(re^{it})}{\partial t} ={\rm Re}\bigg(\frac{z F_z(z)- \bar{z} F_{\bar z}(z)}{F(z)}\bigg) > 0
\end{eqnarray*}
for all $z = r\, e^{it}$ and $r\in (0, 1)$.
\edefe

For a complex-valued function $f$ in $D$, its Jacobian $J_f(z)$ is given by $J_f(z)=|f_z(z)|^2-|f_{\overline{z}}(z)|^2$.
We say that a harmonic mapping $f$ is locally univalent and sense-preserving if and only if its Jacobian $J_f(z)>0$ for $z\in D$ (cf. \cite{L1936}). 
For continuously differentiable function $f$, let 
\begin{eqnarray*}
\Lambda_{f}(z)
=|f_{z}(z)|+|f_{\overline{z}}(z)| ~\mbox{ and }~ \lambda_{f}(z)
=\big ||f_{z}(z)|-|f_{\overline{z}}(z)|\big |.
\end{eqnarray*}

Throughout $\ID_{r}=\{z\in\mathbb{C}:\,|z|<r\}$  denotes the open  disk about the origin so that  $\ID:= \ID_1$ is the unit disk.
For the convenience of the reader, let us fix some basic notations.
\bei
\item ${\rm Hol}(\ID)=\{f:\, \mbox{$f$ is analytic in $\mathbb{D}$}$\}
\item $\mathcal{B}_M=\{f \in {\rm Hol}(\ID):\, |f(z)|\leq M ~\mbox{ in $\ID$}\}$
\item $\mathcal{A}_0=\{f \in {\rm Hol}(\ID):\, f(0) =0\}$ and  $\mathcal{A}_1=\{f \in {\rm Hol}(\ID):\, f'(0) =1\}$
\item $\mathcal{A}=\{f \in {\rm Hol}(\ID):\, f(0) =0=f'(0)-1\}:= \mathcal{A}_1 \cap \mathcal{A}_0$
\item $\mathcal{H}=\{f:\, \mbox{$f$ is harmonic in $\mathbb{D}$}\}$
\item $\mathcal{H}_0=\{f\in \mathcal{H}:\, f(0) =0\}$
\item $\mathcal{BH}_M=\{f\in \mathcal{H}:\, |f(z)|\leq M ~\mbox{ in $\ID$}\}$
\item $\mathcal{BH}_M^0= \mathcal{BH}_M\cap \mathcal{H}_0$
\item ${\rm Bi}\mathcal{H}=\{f:\, \mbox{$f$ is biharmonic in $\mathbb{D}$}\}$
\item ${\rm Bi}\mathcal{A}_0=\{f:\, \mbox{$f$ is bi-analytic in $\mathbb{D}$ with $f(0)=0$}\}.$
\eei

\bdefe
{\it
A function $f$ in a family is said to belong to ${\mathcal S}(r;R)$ if it is univalent in $\ID_r$ and the range
$f(\mathbb{D}_r)$ contains a univalent disk  $\mathbb{D}_R$.
}
\edefe

\subsection{Landau and Bloch theorems}
The classical theorem of Landau states that if $f\in \mathcal{B}_M \cap \mathcal{A}$ for some $M>1$, then $f\in {\mathcal S}(r;R)$ with  $r=1/(M+\sqrt{M^2-1})$ and $R=Mr^2$.
This result is sharp, with the extremal function
$f_0(z)=M z \frac{1-Mz}{M-z}$.

 The Bloch theorem asserts the existence of a positive constant number $b$ such that if $f\in \mathcal{A}_1$, then $f(\ID)$ contains a schlicht disk of radius $b$,
that is, a disk of radius $b$ which is the univalent image of some subregion of the unit disk $\ID$. The supremum of all such constants $b$ is called the Bloch constant (see \cite{CGH2000,GK2003}).

In 2000, under a suitable restriction, Chen et al. \cite{CGH2000} first established two non-sharp versions of Landau-type theorems for bounded harmonic mapping on the unit disk which we now recall with the help of our notation.

\begin{Thm}\label{Theo-A} {\rm (\cite[Theorem 3]{CGH2000})}  
If $f\in \mathcal{BH}_M^0$ with the normalization $f_{\overline{z}}(0)=0$ and $f_z(0)=1$, then $f\in {\mathcal S}(r_1;r_1/2)$, where
$$r_1=\frac{\pi^2}{16 m M}\approx\frac{1}{11.105M},
$$
where $m\approx 6.85$ is the minimum of the function $(3-r^2)/(r(1-r^2))$ for $0<r<1$.
\end{Thm}

\begin{Thm}\label{Theo-B} {\rm (\cite[Theorem 4]{CGH2000})}
If $f\in \mathcal{H}_0$ such that $\lambda_f(0)=1$, and $\Lambda_f(z)\leq\Lambda$ for $z\in \ID$, then  $f\in {\mathcal S}(r_2;r_2/2)$, where
$r_2=\frac{\pi}{4(1+\Lambda )}.$
\end{Thm}

Theorems~A and B 
are not sharp. Better estimates were given in \cite{DN2004} and this topic was later dealt by a number of authors (cf. \cite{cg,cpw2011C,CPR2014,CPW2011B,G2006,Huang2014,L2009S,LC2018}). In 2008, Abdulhadi and Muhanna established two versions of Landau-type theorems for certain bounded biharmonic mappings in \cite{AM2008}.
For later developments on this topic, we refer to \cite{CPW2009, cpw2011C, L2008C,LLZ2011,LXY2017,ZL2013}. In particular, sharp versions of Theorem~B 
have been established in \cite{Huang2014,L2009S,LC2018}, and the corresponding sharp versions of Landau-type theorems for normalized bounded biharmonic mappings  have also been established in \cite{LL2019}.

\begin{Thm}\label{Theo-C} {\rm (\cite[Theorem 3.1]{LL2019})}\quad  Suppose that $\Lambda_1\geq0$ and $\Lambda_2>1$. Let $F\in {\rm Bi}\mathcal{H}$
and  $F(z)=|z|^2G(z)+H(z)$, where
$G,H\in \mathcal{H}_0$, $\lambda_{F}(0)=1,~\Lambda_{G}(z)\leq\Lambda_1$ and $\Lambda_{H}(z)<\Lambda_2$ for all $z\in \ID$.
Then $F\in{\mathcal S}(r_3;R_3)$, where $r_3$ is the unique root in $(0,\,1)$ of the equation
\begin{eqnarray*}
\Lambda_2\, \frac{1-\Lambda_2 r}{\Lambda_2-r}-3\Lambda_1 r^2=0,
\end{eqnarray*}
and
\begin{eqnarray*}
R_3=\Lambda_2^2r_3+\left(\Lambda_2^3-\Lambda_2\right)\ln\left(1-\frac{r_3}{\Lambda_2}\right)-\Lambda_1 r_3^3.
\end{eqnarray*}
This result is sharp.
\end{Thm}
\begin{Thm}\label{Theo-D} {\rm (\cite[Theorem 3.3]{LL2019})} Suppose that
 $\Lambda\geq 0$. Let $F\in {\rm Bi}\mathcal{H}$ and  $F(z)=|z|^2G(z)+H(z)$, where
$G,H\in \mathcal{H}_0$, $\lambda_{F}(0)=1, \Lambda_{G}(z)\leq\Lambda$ and $\Lambda_H(z)\leq 1$ for all $z\in \ID$. Then
$F\in{\mathcal S}(r_4;R_4)$, where
$$
r_4= \left \{ \begin{array}{rl}
1 &\mbox{ when $\ds 0\leq\Lambda\leq\frac{1}{3}$}\\
\ds \frac{1}{\sqrt{3\Lambda}} &\mbox{ when $\ds \Lambda>\frac{1}{3}$},
\end{array}
\right .
$$
and $R_4=r_4-\Lambda r_4^{3}$. This result is sharp.
\end{Thm}

However, the sharp version of Landau-type theorem for normalized bounded harmonic mappings or Theorem~A 
for the case of the bound $M>1$  has not been established. In 2022, Abdulhadi and Hajj established
the following non-sharp Landau-type theorem for certain bounded bi-analytic functions.

\begin{Thm}\label{Theo-E}\cite{AH2022}  Let $F\in {\rm Bi}\mathcal{A}_0$ and
$F(z)=\bar{z}G(z)+H(z)$, where $G, H\in \mathcal{A}\cap \mathcal{B}_M $ for some $M>0$. Then, $F\in{\mathcal S}(r_5;R_5)$,
where
$$ r_5=1-\sqrt{\frac{2M}{2M+1}}  ~\mbox{ and }~ R_5=r_5-r_5^2-M\frac{r_5^2+r_5^3}{1-r_5} .
$$
\end{Thm}

Theorem~E 
is not sharp too.

\subsection{Two natural question on Landau-type theorem} From the discussion above, a couple of natural questions arise.

\bprob\label{HLP-prob1}
Can we establish some sharp versions of Landau-type theorems for certain bounded bi-analytic functions?
\eprob

\bprob\label{HLP-prob2}
Whether we can further generalize and/or improve Theorem~E? 
\eprob


The paper is organized as follows. In Section \ref{HLP-sec2}, we present statements of four theorems out of which one of them improves Theorem~E. 
In addition, we provide several sharp versions of Landau-type theorems for certain bounded bi-analytic functions, which provide an affirmative answer to Problems \ref{HLP-prob1} and \ref{HLP-prob2}. 
In particular, as consequence, we also obtain four sharp versions of Landau-type theorems for certain bounded biharmonic mappings. In Section \ref{HLP-sec3}, we state  a couple of lemmas which are needed for the proofs of
main results in Section \ref{HLP-sec4}.

\section{Statement of Main Results and Remarks}\label{HLP-sec2}

We first establish the following sharp version of Landau-type theorem for certain subclass of bounded bi-analytic functions.
\bthm\label{HLP-th1}
Suppose that $\Lambda_1\geq0$ and $\Lambda_2>1$. Let $F\in {\rm Bi}\mathcal{A}_0$ and $F(z)=\bar{z}G(z)+H(z)$, where $G\in \mathcal{A}_0$, $H\in \mathcal{A}$, $|G'(z)|\leq\Lambda_1$ and
$|H'(z)|<\Lambda_2$ for all $z\in \ID$. Then $F\in{\mathcal S}(\rho_1;\sigma_1)$, where
\begin{eqnarray}
\rho_1=\frac{2\Lambda_2}{\Lambda_2(2\Lambda_1+\Lambda_2)+\sqrt{\Lambda_2^2(2\Lambda_1+\Lambda_2)^2-8\Lambda_1\Lambda_2}},
\label{liu21}
\end{eqnarray}
and
\begin{eqnarray}
\sigma_1=F_1(\rho_1), \quad F_1(z) = \Lambda_2^2 z-\Lambda_1|z|^2+\left(\Lambda_2^3-\Lambda_2\right)\ln\bigg(1-\frac{z}{\Lambda_2}\bigg).
\label{liu22}
\end{eqnarray}
This result is sharp, with an extremal function given by $F_1(z)$.
\ethm


For the case $\Lambda_1\geq0$ and $\Lambda_2=1$,  we will prove the following sharp version of Landau-type theorem for certain subclass of bounded bi-analytic functions.

\bthm \label{HLP-th2}
Suppose that $\Lambda\geq 0$. Let $F\in {\rm Bi}\mathcal{A}_0$ and $F(z)=\bar{z}G(z)+H(z)$, where $G\in \mathcal{A}_0$, $H\in \mathcal{A}$,
$|G'(z)|\leq\Lambda$, and $|H(z)|< 1$ or $|H'(z)|\leq 1$ for all $z\in \ID$. Then $F\in{\mathcal S}(\rho_2;\sigma_2)$, where
$$
\rho_2= \left \{ \begin{array}{rl}
1 &\mbox{ when $\ds 0\leq\Lambda\leq\frac{1}{2}$}\\
\ds \frac{1}{2\Lambda} &\mbox{ when $\ds \Lambda>\frac{1}{2}$},
\end{array}
\right .
$$
and $\sigma_2=\rho_2-\Lambda \rho_2^2$. This result is sharp.
\ethm

\br\label{HLP-re2}
Note that $G\in\mathcal{A}_0$ implies that $G(z)=z G_1(z)$ with $G_1(z)$ being analytic in $\ID$. Thus the bi-analytic function $F(z)=\bar{z}G(z)+H(z)$  reduces to the form
$F(z)=|z|^2G_1(z)+H(z)$ which is clearly a biharmonic mappings. Hence, we conclude the following corollaries from Theorems \ref{HLP-th1} and \ref{HLP-th2}. 
\er

\bcor\label{HLP-cor1}
Suppose that $\Lambda_1\geq0$ and $\Lambda_2>1$. Let $F(z)=|z|^2G(z)+H(z)$ belong to ${\rm Bi}\mathcal{H}$, where $G\in {\rm Hol}(\ID)$ and $H\in \mathcal{A}$.

\begin{enumerate}
\item If $|G(z)+zG'(z)|\leq\Lambda_1$, and $|H'(z)|<\Lambda_2$ for all $z\in \ID$, then $F\in{\mathcal S}(\rho_1;\sigma_1)$ where $\rho_1$  and $\sigma_1$ are
given by \eqref{liu21} and \eqref{liu22}, respectively. This result is sharp, with an extremal function  $F_1(z)$ given by \eqref{liu22}.

\item If $|G(z)+zG'(z)|\leq\Lambda_1$, and $|H(z)|< 1$ or $|H'(z)|\leq 1$ for all $z\in \ID$, then $F\in{\mathcal S}(\rho_2;\sigma_2)$ where $\rho_2$  and $\sigma_2$ are as in Theorem \ref{HLP-th2}.
This result is sharp, with an extremal function given by $F_2(z)=\Lambda _1 |z|^2+z$.
\end{enumerate}
\ecor


If we replace the condition ``$|G(z)+zG'(z)|\leq\Lambda_1$ for all $z\in \ID$'' by the conditions ``$G(0)=0$ and $|G'(z)|\leq\Lambda_1$ for all $z\in \ID$" in Corollary \ref{HLP-cor1}, then, by Theorems~C and D, 
we have the following sharp versions of Landau-type theorems for the special subclasses of bounded biharmonic mappings.

%

\bcor\label{HLP-cor3}
Suppose that $\Lambda_1\geq0$ and $\Lambda_2>1$.  Let $F(z)=|z|^2G(z)+H(z)$ belong to ${\rm Bi}\mathcal{H}$, where $G\in {\rm Hol}(\ID)$ and $H\in \mathcal{A}$.
\begin{enumerate}
\item  If $|G'(z)|\leq\Lambda_1$, and $|H'(z)|<\Lambda_2$ for all $z\in \ID$, then $F\in{\mathcal S}(r_3;R_3)$, where $r_3$ and $R_3$ are as in Theorem~C. 
This result is sharp, with an extremal function given by 
\begin{eqnarray*}
F_0(z)
&=&\Lambda_2^2 z-\Lambda_1|z|^2z+\left(\Lambda_2^3-\Lambda_2\right)\ln\bigg(1-\frac{z}{\Lambda_2}\bigg) .
\end{eqnarray*}
\item If $|G'(z)|\leq\Lambda_1$, and $|H(z)|< 1$ or $|H'(z)|\leq 1$ for all $z\in \ID$, then $F\in{\mathcal S}(r_4;R_4)$, where $r_4$ and $R_4$ are as in Theorem~D. 
This result is sharp. This result is sharp, with an extremal function given by $F_2(z)=\Lambda _1 |z|^2+z$.
\end{enumerate}

\ecor

Now we improve Theorem~A 
by establishing the following results.

\bthm\label{HLP-th3}
Let $F\in {\rm Bi}\mathcal{A}_0$ and $F(z)=\bar{z}G(z)+H(z)$, where $G\in \mathcal{B}_{M_1}\cap  \mathcal{A}$ and $H\in \mathcal{B}_{M_2}\cap  \mathcal{A}$ for some $M_1>0$ and $M_2>0$.
Then $F\in{\mathcal S}(\rho_3;\sigma_3)$, where $\rho_3$ is the unique root in $(0,1)$ of the equation
\begin{eqnarray}
1-\Big(M_2-\frac{1}{M_2}\Big)\frac{2r-r^2}{(1-r)^2}-\Big(M_1-\frac{1}{M_1}\Big)\frac{(3-2r)r^2}{(1-r)^2}-2r=0,
\label{liu24}
\end{eqnarray}
and
$$
\sigma_3=\rho_3-\rho_3^2-\Big(M_2-\frac{1}{M_2}\Big)\frac{\rho_3^2}{1-\rho_3}-\Big(M_1-\frac{1}{M_1}\Big)\frac{\rho_3^3}{1-\rho_3}.
$$
\ethm

\br
If we set $M_1=M_2=1$ in Theorem \ref{HLP-th3}, then it is clear that $G(z)=z$ and $H(z)=z$ by Schwarz lemma. Thus, $\rho_3=\frac{1}{2}$ and $\sigma_3=\frac{1}{4}$ are sharp, with an extremal function  $F_3(z)=|z|^2+z$.
Moreover, if we set $M_1=M_2=M$ in Theorem \ref{HLP-th3}, then one can easily gets an improved version of Theorem~E. 
\er

Furthermore, as with Remark \ref{HLP-re2}, we easily have the following.

\bcor\label{HLP-cor6}
Let $F(z)=|z|^2G(z)+H(z)$ belong to ${\rm Bi}\mathcal{H}$, where $G\in \mathcal{B}_{M_1}\cap  \mathcal{A}$ and $H\in \mathcal{B}_{M_2}\cap  \mathcal{A}$ for some $M_1>0$ and $M_2>0$.
Then $F\in{\mathcal S}(\rho_3;\sigma_3)$,  where $\rho_3$   and $\sigma_3$ are as in Theorem \ref{HLP-th3}.
\ecor

\br
Again, if $M_1=M_2=1$,  then we have $\rho_3=\frac{1}{2}$ and $\sigma_3=\frac{1}{4}$ with an extremal function $F_3(z)=|z|^2+z$.
\er

Finally, we improve Theorem \ref{HLP-th3} by establishing the following theorem.

\bthm\label{HLP-th4}
Let $F\in {\rm Bi}\mathcal{A}_0$ and $F(z)=\bar{z}G(z)+H(z)$, where $0\not\equiv G  \in \mathcal{B}_{M_1}\cap  \mathcal{A}$ and $H\in \mathcal{B}_{M_2}\cap  \mathcal{A}$ for some $M_1>0$ and $M_2>0$.
Then $F$ is sense-preserving, univalent and fully starlike in the disk $\ID_{\rho_3}$, where $\rho_3$ is the unique root in $(0,\,1)$ of Eq. \eqref{liu24}.
\ethm

\section{Key lemmas}\label{HLP-sec3}

In order to prove our main results, we need the following lemmas which play a key role in establishing the subsequent results in Section \ref{HLP-sec4}.

\blem\label{HLP-lem23}
Let $H\in \mathcal{A}_{1}$ and $|H'(z)|<\Lambda$ for all $z\in \ID$ and for some $\Lambda>1$.
\begin{enumerate}
\item For all $z_1,z_2\in \ID_r\, (0<r<1, z_1\neq z_2)$, we have
$$
|H(z_1)-H(z_2)|=\bigg|\int_{\gamma }H'(z)\,dz\bigg|\geq\Lambda\, \frac{1-\Lambda r}{\Lambda-r}\, |z_1-z_2|,
$$
where $\gamma =[z_1,z_2]$ denotes the closed line segment joining $z_1$ and $z_2$.

\item  For $z'\in \partial \ID_r\, (0<r<1)$ with $w'=H(z')\in H(\partial \ID_r)$ and $|w'|=\min\left\{|w|:\,w\in H\left(\partial \ID_r\right)\right\}$, set $\gamma _0=H^{-1}(\Gamma _0)$ and
$\Gamma _0= [0,w'] $ denotes the closed line segment joining the origin and $w'$. Then we have
$$
|H(z')| 
\geq\Lambda\int_{0}^{r}\frac{\frac{1}{\Lambda}-t}{1-\frac{t}{\Lambda}}\,dt=\Lambda^2 r+(\Lambda^3-\Lambda)\ln\left(1-\frac{r}{\Lambda}\right).
$$
\end{enumerate}
\elem

\bpf Set $\omega(z)=H'(z)/\Lambda$, $z\in \ID$. Then $\omega \in \mathcal{B}_1$ with $\alpha:=\omega(0)=\frac{H'(0)}{\Lambda}=\frac{1}{\Lambda}.$
Using Schwarz-Pick Lemma, we have
$$
\frac{\frac{1}{\Lambda}-r}{1-\frac{r}{\Lambda}}= \frac{\alpha -r}{1-\alpha r} \leq {\rm Re}\, \omega(z) \leq |\omega(z)|\leq  \frac{\alpha +r }{1+\alpha r} , \quad z\in \ID_r.
$$

(1) Fix $z_1,z_2\in \ID_r\, (0<r<1)$ with $z_1\neq z_2$, set $\theta_0=\arg(z_2-z_1)$. Then 
\begin{eqnarray*}
\hspace{1cm}|H(z_1)-H(z_2)|&=&\bigg|\int_{\overline{\gamma }}H'(z)\,dz\bigg|=\left| \int_{\gamma}\Lambda\, \omega(z)e^{i\theta_0}\,|dz|\right|\\
&\geq& \Lambda\int_{\gamma}{\rm Re}\, \omega(z)|dz|\\
&\geq&\Lambda\, \int_{\gamma}\frac{\frac{1}{\Lambda}-r}{1-\frac{r}{\Lambda}}\,|dz|=\Lambda\, \frac{1-\Lambda r}{\Lambda-r}\, |z_1-z_2|.
\end{eqnarray*}

(2) For $z'\in \partial \ID_r\, (0<r<1)$ with $w'=H(z')\in H(\partial \ID_r)$,    $|w'|=\min\left\{|w|:\,w\in F\left(\partial \ID_r\right)\right\}$ and $\Gamma _0=[0,w]$,
set $\gamma _0=H^{-1}(\Gamma _0)$ so that
\begin{eqnarray*}
\hspace{1.5cm}|H(z')|&=&|w'|=
\int_{\gamma_0}|H'(\zeta)|\,|d\zeta|=\Lambda\int_{\gamma_0}|\omega (\zeta)|\,|d\zeta|\\
&\geq&\Lambda\int_{0}^{r} \frac{\frac{1}{\Lambda}-t}{1-\frac{t}{\Lambda}}dt=\Lambda^2 r+(\Lambda^3-\Lambda)\ln\left(1-\frac{r}{\Lambda}\right) 
\end{eqnarray*}
and the proof is complete.
\epf

\blem\label{HLP-lem21}{\rm (Carlson lemma, \cite{Ca1940})}\quad
If $F\in \mathcal{B}_{1}$ and $F(z)=\sum_{n=0}^{\infty}a_nz^{n}$, then the following inequalities hold:

\begin{enumerate}
\item[{\rm (a)}] $|a_{2n+1}|\leq 1-|a_0|^2-\cdots - |a_n|^2,\, n=0, 1, \ldots$.

\item[{\rm (b)}] $|a_{2n}|\leq 1-|a_0|^2-\cdots - |a_{n-1}|^2-\frac{|a_n|^2}{1+|a_0|},\, n=1, 2, \ldots$.
\end{enumerate}
These inequalities are sharp.
\elem


\blem\label{HLP-lem22}
If $f\in \mathcal{B}_M\cap \mathcal{A}_0$ for some $M>0$ and $f(z)=\sum_{n=1}^{\infty}a_nz^n$, then
\begin{enumerate}
\item[{\rm (a)}] $\ds |a_{2n}|\leq M\left [1- \left (\frac{|a_1|^2+ \cdots + |a_n|^2}{M^2}\right )\right ],\, n=1, 2, \ldots$.

\item[{\rm (b)}] $\ds |a_{2n+1 }|\leq M\left [1 - \left (\frac{|a_1|^2 + \cdots + |a_n|^2}{M^2}\right ) - \frac{|a_{n+1}|^2}{M(M+|a_1|)} \right ], \, n=1, 2, \ldots $
\end{enumerate}
In particular if $|a_1|=1$, i.e., if $f\in   \mathcal{B}_M\cap \mathcal{A}$, then $M\geq 1$  and
\begin{equation*}
|a_n|\leq M-\frac{1}{M} ~\mbox{ for  $n=2, 3, \ldots .$} 
\end{equation*}
These inequalities are sharp, with the extremal functions $f_n(z)$, where
\begin{eqnarray*}
f_1(z)=z,\quad f_n(z)=Mz\frac{1-M z^{n-1}}{M-z^{n-1}}=z-\Big(M-\frac{1}{M}\Big)z^n-\sum\limits_{k=3}^\infty\frac{M^2-1}{M^{k-1}}z^{(n-1)(k-1)+1}
\end{eqnarray*}
for $n=2,3,\ldots$.
\elem

\bpf
Setting $g(z)=\frac{f(z)}{M z}$ for $z\in \ID\backslash\{0\}$, and $g(0)=\frac{a_1}{M}$, shows that $g\in \mathcal{B}_1$ and
$$g(z)=  \sum_{n=0}^{\infty}b_nz^{n},
$$
where $b_n= a_{n+1}/M$ for $n\geq 0$. Note that $b_0=a_1/M$. Applying Lemma \ref{HLP-lem21} to the coefficients $b_n$ of $g$ gives the desired inequality.

In particular if $|a_1|=1$,  then we have  $M\geq 1$ and it follows from (a) and (b) that
$$
|a_{n}|\leq M\left (1- \frac{|a_1|^2}{M^2}\right ) =M-\frac{1}{M} ~\mbox{ for $n\geq 2$}~
$$
and it is evident that equalities hold for all $n=2,3,\ldots$ for the functions
$$
f_n(z)=Mz\frac{1-M z^{n-1}}{M-z^{n-1}}=z-\Big(M-\frac{1}{M}\Big)z^n-\sum\limits_{k=3}^\infty\frac{M^2-1}{M^{k-1}}z^{(n-1)(k-1)+1},
$$
and the proof is complete.
\epf 

\blem\label{HLP-lem24}
\quad {\it Let $F(z)=\bar{z}G(z)+H(z)$ be a bi-analytic function of the unit disk $\ID$, where $G(z)=\sum\limits_{n=1}^\infty a_n z^n\not\equiv 0$ and $H(z)=z+\sum\limits_{n=2}^\infty b_n z^n$ are analytic in $\ID$, and satisfy the condition
\begin{eqnarray}
\sum\limits_{n=2}^\infty n|b_n|r^{n-1}+\sum\limits_{n=1}^\infty (n+1)|a_n|r^{n}\leq 1,
\label{liu34}
\end{eqnarray}
for some $r\in (0, 1)$. Then $F(z)$ is sense-preserving, univalent and fully starlike in the disk $\ID_r$.}
\elem

\bpf We may use arguments similar to those in the proof of \cite[Lemma 1]{LP2021}. For the sake of readability,
we provide the details. Elementary computation gives
\begin{eqnarray}
z F_z(z)-\bar{z}F_{\bar{z}}(z)-F(z) 
&=&\bar{z}\sum\limits_{n=1}^\infty (n-2)a_n z^n+\sum\limits_{n=2}^\infty (n-1)b_n z^n.
\label{liu36}
\end{eqnarray}
Evidently, $J_{F}(0)=1$. Now, we fix $r\in (0, 1]$ and find that
\begin{eqnarray*}
|F_z(z)|-|F_{\bar{z}}(z)|&=&\Big|\bar{z}\sum\limits_{n=1}^\infty n a_n z^{n-1}+1+\sum\limits_{n=2}^\infty n b_n z^{n-1}\Big|
-\Big|\sum\limits_{n=1}^\infty a_n z^n\Big|\\
&>& 1-\sum\limits_{n=2}^\infty n|b_n|r^{n-1}-\sum\limits_{n=1}^\infty (n+1)|a_n|r^{n}\geq 0,
\end{eqnarray*}
and therefore,
$J_{F}(z)= (|F_z(z)|+|F_{\bar{z}}(z)|)(|F_z(z)|-|F_{\bar{z}}(z)|)>0$  for $|z|<r$.

Thus, $F$  is sense-preserving in $\ID_r$. Finally, fix $r_0\in (0, r]$ and consider the circle $\partial \ID_{r_0}=\{z:\,|z|=r_0\}$. For $z\in \partial \ID_{r_0}$, it follows from $G(z)=\sum\limits_{n=1}^\infty a_n z^n\not\equiv 0$, (\ref{liu34}) and (\ref{liu36}) that
\begin{eqnarray*}
|z F_z(z)-\bar{z}F_{\bar{z}}(z)-F(z)|&\leq &\sum\limits_{n=1}^\infty |n-2|\,|a_n|\,|z|^{n+1}+\sum\limits_{n=2}^\infty (n-1)|b_n|\,|z|^n\\
&=&|z| \Big(\sum\limits_{n=2}^\infty n|b_n|\,|z|^{n-1}+\sum\limits_{n=1}^\infty (n+1)|a_n|\,|z|^{n}\Big)\\
&&-|a_1|\,|z|^2-3\sum\limits_{n=2}^\infty |a_n|\,|z|^{n+1}-\sum\limits_{n=2}^\infty |b_n|\,|z|^{n}\\
&< & |z|-\sum\limits_{n=2}^\infty |b_n|\,|z|^{n}-|\bar{z}|\sum\limits_{n=1}^\infty |a_n|\,|z|^{n}\\
&\leq & |H(z)|-|\bar{z}G(z)|\leq |F(z)|,
\end{eqnarray*}
which implies that
$$
\bigg|\frac{z F_z(z)- \bar{z} F_{\bar z}(z)}{F(z)}-1\bigg|<1\quad\mbox{ for } |z|=r_0.
$$
Thus, we obtain that $F$ is univalent on $\partial \ID_{r_0}$, and it maps $\partial \ID_{r_0}$ onto a starlike curve. Hence, by the sense-preserving property and the degree principle, we see that $F$ is univalent in $\ID_{r_0}$. Since $r_0\in (0, r]$ is arbitrary, we conclude that $F$ is univalent and fully starlike in $\ID_r$. The proof is complete. \epf

\section{Proofs of the main results}\label{HLP-sec4}

\subsection{Proof of Theorem \ref{HLP-th1}}
By the assumption on $G\in \mathcal{A}_0$, we have
\begin{eqnarray}
|G(z)|=\bigg|\int_{[0,z]}G'(z)\, dz\bigg|\leq \int_{[0,z]}|G'(z)|\,|dz|\leq \Lambda_1 |z|,\ \ z\in \ID.
\label{liu41}
\end{eqnarray}

We first prove that $F$ is univalent in the disk $\ID_{\rho_1}$. Choose, for all $z_1,z_2\in \ID_r\, (0<r<\rho_1$, $z_1\neq z_2)$, where $\rho_1$ is defined by \eqref{liu21}. As  $H'(0)=1$, $|G'(z)|\leq\Lambda_1$ and $|H'(z)|<\Lambda_2$ for all $z\in \ID$, we obtain from Lemma \ref{HLP-lem23} that
\begin{eqnarray*}
|F(z_2)-F(z_1)|&=&\left|\int_{[z_1,z_2]}F_z(z)\,dz+F_{\bar{z}}(z)\,d\bar{z}\right|
\, =\left|\int_{[z_1,z_2]}\left(\bar{z}G'(z)+H'(z)\right)dz+G(z)\,d\bar{z}\right|\nonumber\\
&\geq&\bigg|\int_{[z_1,z_2]}H'(z)\,dz\bigg|-\bigg|\int_{[z_1,z_2]}\bar{z}G'(z)\,dz+G(z)\,d\bar{z}\bigg|\\
&\geq&
|z_1-z_2|\left(\Lambda_2 \frac{1-\Lambda_2 r}{\Lambda_2-r}-2\Lambda_1r\right)\nonumber\\
&=&|z_1-z_2|\cdot\frac{2\Lambda_1r^2-\Lambda_2(2\Lambda_1+\Lambda_2)r+\Lambda_2}{\Lambda_2-r}\nonumber\\
&=&|z_1-z_2|\frac{2\Lambda_1(r-\rho_1)(r-A)}{\Lambda_2-r}\nonumber
\end{eqnarray*}
which is positive if $r<\rho_1$, where
\begin{eqnarray}
A=\frac{\Lambda_2(2\Lambda_1+\Lambda_2)+\sqrt{\Lambda_2^2(2\Lambda_1+\Lambda_2)^2-8\Lambda_1\Lambda_2}}{4\Lambda_1}.\nonumber
\end{eqnarray}

This  proves the univalency of $F$ in the disk $\ID_{\rho_1}$.

Next, we prove that $F(\ID_{\rho_1}) \supseteq \ID_{\sigma_1}$, where $\sigma_1$ is defined by \eqref{liu22}. First we note that $F(0)=0$, for $z'\in \partial \ID_{\rho_1}$ with $w'=F(z')\in F(\partial \ID_{\rho_1})$ and $|w'|=\min\left\{|w|:\,w\in F\left(\partial \ID_{\rho_1}\right)\right\}$. By \eqref{liu41} and Lemma \ref{HLP-lem23}, we have that
\begin{eqnarray*}
|w'|&=&\big|\bar{z'}G(z')+H(z')\big|\geq |H(z')|-\Lambda_1 \rho_1^2
\geq h_0(\rho_1)=\sigma_1,
\end{eqnarray*}
which implies that $F(\ID_{\rho_1})\supseteq \ID_{\sigma_1}$, where
\begin{eqnarray}
h_0(x)=\Lambda_2^2x-\Lambda_1x^2+\left(\Lambda_2^3-\Lambda_2\right)\ln\left(1-\frac{x}{\Lambda_2}\right), \quad  x\in [0, 1].
\label{liu43}
\end{eqnarray}

Now, we prove the sharpness of $\rho_1$ and $\sigma_1$. To this end, we consider the bi-analytic function $F_1(z)$ which is given by (\ref{liu22}). It is easy to verify that $F_1(z)$ satisfies the hypothesis of Theorem \ref{HLP-th1}, and thus, we have that $F_1(z)$ is univalent in $\ID_{\rho_1}$, and $F_1(\ID_{\rho_1}) \supseteq \ID_{\sigma_1}$.

To show that the radius $\rho_1$ is sharp, we need to prove that $F_1(z)$ is not univalent in $\ID_r$ for each $r\in (\rho_1, 1]$. In fact for the real differentiable function $h_0(x)$
given above, we have
$$
h_0'(x)=\frac{2\Lambda_1x^2-\Lambda_2(2\Lambda_1+\Lambda_2)x+\Lambda_2}{\Lambda_2-x}
$$
which is continuous and strictly decreasing on $[0, 1]$ with $h_0'(\rho_1)=0$. It follows that $h_0'(x)=0$ for $x\in [0, 1]$ if and only if $x=\rho_1$. So $h_0(x)$ is strictly increasing on $[0, \rho_1)$
and strictly decreasing on $[\rho_1, 1]$. Since $h_0(0)=0$, there is a unique real $\rho_1'\in (\rho_1, 1]$ such that $h_0(\rho_1')=0$ if $h_0(1)\leq 0$, and
\begin{eqnarray}
\sigma_1=\Lambda_2^2\rho_1+\left(\Lambda_2^3-\Lambda_2\right)\ln\left(1-\frac{\rho_1}{\Lambda_2}\right)-\Lambda_1\rho_1^2=h_0(\rho_1)>h_0(0)=0.
\label{liu44}
\end{eqnarray}

For every fixed $r\in (\rho_1, 1]$, set $x_1=\rho_1+\varepsilon$, where
\begin{equation*}
\varepsilon=\left\{
\begin{array}{lll}
\ds\min\left\{\frac{r-\rho_1}{2}, \frac{\rho_1'-\rho_1}{2}\right\} && \mbox{ if } h_0(1)\leq 0,\\
\ds \frac{r-\rho_1}{2} && \mbox{ if } h_0(1)>0.
\end{array}
\right.
\end{equation*}
By the mean value theorem, there is a unique $\delta\in (0, \rho_1)$ such that $x_2:=\rho_1-\delta\in (0, \rho_1)$ and $h_0(x_1)=h_0(x_2)$.

Let $z_1=x_1$ and $z_2=x_2$. Then $z_1,\ z_2\in \ID_r$ with $z_1\neq z_2$ and observe that
\begin{eqnarray*}
F_1(z_1)=F_1(x_1)=h_0(x_1)=h_0(x_2)=F_1(z_2).
\end{eqnarray*}
Hence $F_1$ is not univalent in the disk $\ID_r$ for each $r\in (\rho_2, 1]$, and thus, the radius $\rho_1$ is sharp.

Finally, note that $F_1(0)=0$ and picking up $z'=\rho_1\in \partial \ID_{\rho_1}$, by (\ref{liu22}), (\ref{liu43}) and (\ref{liu44}), we have
\begin{eqnarray*}
|F_1(z')-F_1(0)|=|F_1(\rho_1)|=|h_0(\rho_1)|=h_0(\rho_1)=\sigma_1.
\end{eqnarray*}
Hence, the radius $\sigma_1$ of the schlicht disk is also sharp.\hfill $\Box$

\subsection{Proof of Theorem \ref{HLP-th2}}
The assumption on $H$, namely, $H\in {\mathcal B}_1\cap {\mathcal A}$,   clearly gives that $H(z)\equiv z$ in $\ID$ (by Schwarz's lemma).
Thus, $F$ reduces to the form $F(z)=\bar{z} G(z)+z$.

Now we prove $F$ is univalent in the disk $\ID_{\rho_1}$. To this end, for any $z_1,z_2\in \ID_r\, (0<r<\rho_2)$ with $z_1\neq z_2$, by the condition $G(0)=0$ and $|G'(z)|\leq\Lambda$ for all $z\in \ID$, and (\ref{liu41}), it follows that $|G(z)|\leq \Lambda |z|$ in $\ID$.
Consequently,
\begin{eqnarray*}
|F(z_1)-F(z_2)|&\geq& 
|z_1-z_2| -\bigg|\int_{[z_1,z_2]}\bar{z}G'(z)dz+G(z)d\bar{z}\bigg|\\
&\geq &|z_1-z_2|(1-2\Lambda r)>0
\end{eqnarray*}
which proves the univalency of $F$ in the disk $\ID_{\rho_2}$, where $\rho_2$ is given in the statement of the theorem.

Noticing that $F(0)=0$, for any $z=\rho_2 e^{i\theta}\in \partial \ID_{\rho_2}$, we have
\begin{eqnarray*}
|F(z)|&=&|\bar{z} G(z)+z|\geq |z|-\rho_1|G(z)| \geq \rho_2-\Lambda \rho_2^2=\sigma_2.
\end{eqnarray*}
Hence, $F(\ID_{\rho_2})$ contains a schlicht disk $\ID_{\sigma_2}$.

Finally, for $F_2(z)=\Lambda|z|^2+z$,
a direct computation verifies that $\rho_2$ and $\sigma_2$ are sharp. This completes the proof.\hfill $\Box$

\subsection{Proof of Theorem \ref{HLP-th3}}
As  $G\in \mathcal{B}_{M_1}\cap \mathcal{A}$ and $H\in \mathcal{B}_{M_2}\cap \mathcal{A}$ by assumption, we may write
$$
G(z)=\sum\limits_{n=1}^\infty a_n z^n ~\mbox{ and }~ H(z)=\sum\limits_{n=1}^\infty b_n z^n
$$
where  $a_1=b_1=1$, and it follows from Lemma \ref{HLP-lem22} that
\begin{eqnarray}
|a_n|\leq M_1-\frac{1}{M_1}~\mbox{ and }~ |b_n|\leq M_2-\frac{1}{M_2} ~\mbox{ for all $n\geq 2$}.
\label{liu46}
\end{eqnarray}

We first prove that $F$ is univalent in the disk $\ID_{\rho_3}$, where $\rho_3$ is defined by \eqref{liu24}. Indeed, for all $z_1,z_2\in \ID_r\, (0<r<\rho_3$, $z_1\neq z_2)$, we see that (with $\gamma =[z_1,z_2]$)
\begin{eqnarray*}
&&|F(z_2)-F(z_1)|=\left|\int_{\gamma}F_z(z)dz+F_{\bar{z}}(z)d\bar{z}\right|\nonumber\\
&\geq&\bigg|\int_{\gamma}H'(0)dz\bigg|-\bigg|\int_{\gamma}(H'(z)-H'(0))dz\bigg|
-\bigg|\int_{\gamma}\bar{z}G'(z)dz+G(z)d\bar{z}\bigg|\\
&\geq& |z_1-z_2|\bigg[1-\sum\limits_{n=2}^\infty n|b_n|r^{n-1}-\sum\limits_{n=1}^\infty (n+1)|a_n|r^{n}\bigg]\nonumber\\
&\geq&|z_1-z_2|\bigg [1-\Big(M_2-\frac{1}{M_2}\Big)\sum\limits_{n=2}^\infty nr^{n-1}-\Big(M_1-\frac{1}{M_1}\Big)\sum\limits_{n=2}^\infty (n+1)r^{n}-2r\bigg ]\nonumber\\
&=&|z_1-z_2|\Big [1-\Big(M_2-\frac{1}{M_2}\Big)\frac{2r-r^2}{(1-r)^2}-\Big(M_1-\frac{1}{M_1}\Big)\frac{(3-2r)r^2}{(1-r)^2}-2r\Big ]>0.\nonumber
\end{eqnarray*}
This implies $F(z_1)\neq F(z_2)$, which proves the univalency of $F$ in the disk $\ID_{\rho_3}$.

Next, we prove that $F(\ID_{\rho_3}) \supseteq \ID_{\sigma_3}$, where $\sigma_3$ is as in the statement. Indeed, note that $F(0)=0$ and for any $z'\in \partial \ID_{\rho_3}$ with $w'=F(z')\in F(\partial \ID_{\rho_3})$, 
it follows from (\ref{liu46}) that
\begin{eqnarray*}
|w'|&=&\big|\bar{z'}G(z')+H(z')\big|\geq |H(z')|-\rho_3 |G(z')|\\
&\geq& |z'|-\sum\limits_{n=2}^\infty |b_n|\,|z'|^{n}-\rho_3\sum\limits_{n=1}^\infty |a_n|\,|z'|^{n}\\
&\geq&\rho_3-\rho_3^2-\Big(M_2-\frac{1}{M_2}\Big)\frac{\rho_3^2}{1-\rho_3}-\Big(M_1-\frac{1}{M_1}\Big)\frac{\rho_3^3}{1-\rho_3}=\sigma_3,
\end{eqnarray*}
which implies that $F(\ID_{\rho_3})\supseteq \ID_{\sigma_3}$.
\hfill $\Box$

\subsection{Proof of Theorem \ref{HLP-th4}}
We apply Lemmas \ref{HLP-lem22} and  \ref{HLP-lem24}. Now, by the assumption and the method of proof of Theorem \ref{HLP-th3},
the inequalities in \eqref{liu46} hold and thus, we have
\begin{eqnarray*}
&&\sum\limits_{n=2}^\infty n|b_n|r^{n-1}+\sum\limits_{n=1}^\infty (n+1)|a_n|r^{n}\\
&\leq& \Big(M_2-\frac{1}{M_2}\Big)\sum\limits_{n=2}^\infty nr^{n-1}+\Big(M_1-\frac{1}{M_1}\Big)\sum\limits_{n=2}^\infty (n+1)r^{n}+2r\leq 1
\end{eqnarray*}
for $r\leq \rho_3$. Hence the desired conclusion of Theorem \ref{HLP-th4} follows from Lemma \ref{HLP-lem24}. \hfill $\Box$



\subsection*{Acknowledgments}
This research of the first  author is partly supported by Guangdong Natural Science Foundations (Grant No. 2021A1515010058).

\subsection*{Conflict of Interests}
The authors declare that they have no conflict of interest, regarding the publication of this paper.

\subsection*{Data Availability Statement}
The authors declare that this research is purely theoretical and does not associate with any datas.

\end{document}